\documentclass{amsart}
\usepackage{amsmath}
\usepackage{amsthm}
\usepackage{amssymb}

\usepackage{fullpage}

\usepackage{amsfonts}
\usepackage{latexsym}
\usepackage{mathrsfs}
\usepackage{hyperref}

\numberwithin{equation}{section}
\numberwithin{subsection}{section}

\newtheorem*{thm*}{Theorem}

\theoremstyle{definition}

\newcommand{\rif}{\mathcal{RIF}}

\usepackage{color}
\title{A note on homotopies of rational matrix inner functions}

\author{Michael T. Jury}
\address{Department of Mathematics, University of Florida, PO Box 118105, Gainesville FL 32611-8105 USA}
\email{mjury@ufl.edu}

\begin{document}
\begin{abstract}We show that when $m>n$, the space of $m\times n$-matrix-valued rational inner functions in the disk is path connected. \end{abstract}
\thanks{Partially supported by NSF grant DMS-2154494}
\thanks{This article resulted from a question posed to the author by Prof. Jonathan Adams, Department of Economics, University of Florida, in connection with a problem arising in the analysis of macroeconomic models; see \cite{adams-working}}

\maketitle


A matrix-valued rational function is an $m\times n$ matrix $W(z)$ each of whose entries is a rational function $w_{ij}(z)$ of the complex variable $z$. Thus $W(z)$ is an $m\times n$ matrix valued function defined at all but (at most) finitely many points of the complex plane $\mathbb C$. 

We let $\|W\|_\infty$ denote the supremum of $\|W(z)\|$ over the open unit disk $|z|<1$, here $\|W(z)\|$ is the usual operator norm of the linear transformation $W(z)$ acting between the Euclidean spaces $\mathbb C^n$ and $\mathbb C^m$.  For rational $W$, if $\|W\|_\infty<\infty$, then $W$ extends continuously to the closed disk $|z|\leq 1$, and conversely. (Evidently this occurs if and only if $W$ has no poles in $|z|\leq 1$, we will be working only with such functions.) We say an $m\times n$ rational matrix function is {\em inner} if $\|W\|_\infty\leq 1$ and $W(e^{i\theta})^*W(e^{i\theta})=I_n$ for all $\theta\in [0, 2\pi]$. (Note that this condition forces $m\geq n$.) We will let $\rif(m,n)$ denote the set of all $m\times n$ matrix rational inner functions.  The set $\rif(m,n)$ is equipped with the (metric) topology induced by the norm $\|\cdot\|_\infty$, which it inherits as a subset of the continuous $m\times n$ matrix valued functions in the disk, this coincides with the topology of uniform convergence in the closed disk $|z|\leq 1$. Rational matrix functions (and their inner-outer factorizations) play a fundamental role in many problems of systems theory, automatic control, and prediction theory, among other applications. 
(See for example \cite{MR2363355} and its references.) The purpose of this note is to prove the following:

\begin{thm*}\label{thm:main} If $m>n$ then the metric space $\rif(m,n)$ is path connected.
\end{thm*}

{\bf Remark:} It is easy to see that in the square case, $\rif(m,m)$ is not path connected. Indeed, by considering the winding number of the function $\det W(e^{i\theta})$ about the origin, one sees that, for example, $W(z)=zI_m$ cannot be joined to $I_m$ by a path lying within $\rif(m,m)$. 

\begin{proof}
Since we are assuming $m>n$, it will be helpful to write elements of $\rif(m,n)$ in block form as columns
\[
W(z) = \begin{pmatrix} X(z) \\ Y(z)\end{pmatrix}
\]
where $X(z)$ is an $n\times n$ rational matrix function and $Y(z)$ is $(m-n)\times n$. The fact that $W$ is inner is then expressed by the condition $X(e^{i\theta})^*X(e^{i\theta}) +Y(e^{i\theta})^*Y(e^{i\theta})\equiv I_n$.  

We will prove that every $W=\begin{pmatrix} X\\ Y\end{pmatrix}\in \rif(m,n)$ can be joined to $\begin{pmatrix} I_n \\ O_{(m-n)\times n}\end{pmatrix}$ by a path in $\rif(m,n)$, this evidently proves the theorem. This in turn is accomplished in two steps: first we prove that for any $W\in \rif(m,n)$, there is a square matrix rational inner function $\Phi(z)\in \rif(n,n)$ such that there is a path in $\rif(m,n)$ joining $W$ to $\begin{pmatrix} \Phi \\ O\end{pmatrix}$. (Here $O$ is the $(m-n)\times n$ zero matrix, henceforth we will drop the size subscripts when they are clear from context.) Then we will show that any such $\begin{pmatrix} \Phi \\ O\end{pmatrix}$ can be joined to $\begin{pmatrix} I\\ O\end{pmatrix}$ in $\rif(m,n)$.

       Since $W^*W\equiv I$ on the circle, the matrix $W(e^{i\theta})$ has full rank $n$ for each $\theta\in[0,2\pi)$. In particular, the matrix $W(1)$ has $n$ linearly independent rows, and by continuity this same set of rows is independent in $W(e^{i\theta})$ for $\theta$ in a neighborhood of $0$. Multiplying $W$ on the left by an $m\times m$ permutation matrix, we may arrange that these are the first $n$ rows. Since the unitary group $\mathcal U(m)$ is path connected, and a unitary times a matrix $RIF$ is again a $RIF$, it follows that the new $W$ with permuted rows is connected by a path in $\rif(m,n)$ to the original $W$. So, we may assume $W=\begin{pmatrix} X\\ Y\end{pmatrix}\in \rif(m,n)$ with $X(e^{i\theta})$ having full rank for $\theta$ in a neighborhood of $0$. The rational matrix function $X$ admits an inner-outer factorization $X=\Phi F$, where $\Phi$ is projection-valued on the circle and $F$ is a matrix outer function satisfying $F^*F=X^*X$ on the unit circle; $F$ will be unique if we additionally impose the condition that $F(0)$ be positive definite (which we do). From the theory of matrix inner-outer factorizations, $F$ is also rational \cite[Section 6.8]{MR1435287}. Since $X(1)$ has full rank, it follows that $\Phi(1)$ has full rank $n$, but then by continuity $\text{rank}(\Phi(e^{i\theta}))=\text{trace}(\Phi(e^{i\theta})^*\Phi(e^{i\theta}))$ is constantly equal to $n$. Thus $\Phi\in \rif(n,n)$. We may then write

      \[
      \begin{pmatrix} X\\ Y\end{pmatrix} = \begin{pmatrix} \Phi & 0 \\ 0 & I\end{pmatrix} \begin{pmatrix} F \\ Y\end{pmatrix}.
            \]
            Since $F^*F=X^*X$ on the circle, it follows that $V:=\begin{pmatrix} F\\ Y\end{pmatrix}$ is inner, i.e. belongs to $\rif(m,n)$. If we show that $V$ can be joined to $\begin{pmatrix} I\\ 0\end{pmatrix}$, then (since multiplication by $diag(\Phi, I)$ will carry $\rif(m,n)$ into itself continuously) it will follow that $W$ can be joined to $\begin{pmatrix} \Phi\\ 0\end{pmatrix}$.

                Now, for $0\leq t\leq 1$ the $n\times n$ matrix function $Q_t(e^{i\theta})=I-t^2 Y(e^{i\theta})^*Y(e^{i\theta})$ takes positive semidefinite values on the unit circle (in fact positive definite values when $0\leq t<1$).  Since $Y$ is a rational matrix function, we can choose a polynomial $p$ of minimal degree with the property that $\tilde{Y}(z) := p(z) Y(z)$ is a polynomial matrix function. (That is, $p$ is a common denominator for the entries of $Y$.) Since $Y$ has no poles in $|z|\leq 1$, this minimal degree common denominator will have no zeroes in $|z|\leq 1$, and we may normalize so that $p(0)>0$.  We then consider the nonnegative matrix-valued trigonometric polynomials $Q_t$ given by
\[
  \tilde{Q}_t(e^{i\theta})= \overline{p(e^{i\theta})}p(e^{i\theta})I_n -t^2 \tilde{Y}(e^{i\theta})^*\tilde{Y}(e^{i\theta}).             
\]
By the Fejer-Riesz theorem for matrix valued trigonometric polynomials \cite[Section 6.6]{MR1435287}, there is an outer (analytic) polynomial matrix function $G_t(z)$, with $\deg G_t=\deg Q_t\leq \max(\deg p, \deg Y)$, such that
\[
\overline{p(e^{i\theta})}p(e^{i\theta})I_n -t^2 \tilde{Y}(e^{i\theta})^*\tilde{Y}(e^{i\theta})  = G_t(e^{i\theta})^*G_t(e^{i\theta}).
\]
This $G_t$ will be unique if we impose the requirement that $G_t(0)$ be positive definite. Doing this, in particular we will have $G_0(z)=p(z)I_n$ and $G_1(z)=p(z)F(z)$. Moreover, the outer factor $G_t$ has the following extremal property: if $R$ is any other matrix function, bounded by $1$ in the disk and which satisfies $R^*R\leq \widetilde{Q}_t$ on the circle, then $R(0)^*R(0)\leq G_t(0)^*G_t(0)$ (this follows from the extremal characterization of matrix outer functions \cite[Theorem C, Section 3.10]{MR1435287}). In addition, since all the $G_t$ have full rank and are outer, it follows that $\det G_t(z)$ is nonvanishing in $|z|<1$ for all $0\leq t\leq 1$.

With these facts in hand we can prove that the map $t\to G_t$ is norm continuous on $[0,1]$. We must show that if $t_n\to t$ then $G_{t_n}\to G_t$ uniformly. Since the norms and degrees of the polynomials $G_t$ are uniformly bounded, by compactness there will be a subsequence $G_{t_{n_k}}$ which converges uniformly in $|z|\leq 1$ to some polynomial matrix function $H(z)$. Next we observe that $G_t(0)^*G_t(0)\geq G_1(0)^*G_1(0)$ for all $0\leq t\leq 1$ (this follows from the fact that by definition $G_1^*G_1\leq G_t^*G_t$ on the circle, and the extremal property of outer functions noted above). We thus have $G_t(0)^*G_t(0)\geq G_1(0)^*G_1(0)=|p(0)|^2F(0)^*F(0)$ for all $t$, and since $F(0)$ is positive definite, it follows that $\det H(0)=\lim_k \det G_{t_{n_k}}(0)\neq 0$. Hence, from Hurwitz's theorem we conclude that $\det H(z)=\lim \det G_{t_{n_k}}(z)$ is nonvanishing in $|z|<1$, so (since $H$ is polynomial) $H(z)$ is outer. But by uniform convergence it follows that $H(0)>0$ and $\overline{p(e^{i\theta})}p(e^{i\theta})I_n -t^2 \tilde{Y}(e^{i\theta})^*\tilde{Y}(e^{i\theta})  = H(e^{i\theta})^*H(e^{i\theta})$ for all $\theta$, so by uniqueness we must have $H=G_t$. Thus, for each fixed sequence $t_n\to t$, every subsequence of $G_{t_n}$ has a subsequence converging to $G_t$, so the full sequence converges to $G_t$, and thus $t\to G_t$ is continuous. If we now put $F_t = p^{-1}G_t$, then each $F_t$ is a rational matrix function satisfying
\[
F_t(e^{i\theta})^*F_t(e^{i\theta}) +t^2Y(e^{i\theta})^*Y(e^{i\theta}) \equiv I_n
\]
(with $F_t(0)$ positive definite) for $0\leq t\leq 1$, and the path $t\to F_t$ is continuous. By construction we have $F_0=I_n$ and $F_1=F$. Thus, the columns $\begin{pmatrix} F_t \\ tY\end{pmatrix}$ will belong to $\rif(m,n)$, and form a path joining $\begin{pmatrix} F \\ Y\end{pmatrix}$ to $\begin{pmatrix} I \\ 0\end{pmatrix}$.
Finally, if we put $X_t=\Phi F_t$, then $W_t:=\begin{pmatrix} X_t \\ tY\end{pmatrix}$ is a continuous path in $\rif(m,n)$ joining $W_0=\begin{pmatrix} \Phi \\ 0\end{pmatrix}$ to $W_1= \begin{pmatrix} X \\ Y\end{pmatrix}$ as desired.  

To carry out the second step of the proof, let $\Phi\in \rif(n,n)$. By \cite{MR0114915} $\Phi$ can be factored as a Blaschke-Potapov product
\[
\Phi(z) = U \left(\prod_{k=1}^N \left(b_k(z)P_k +(I-P_k)\right)\right)V
\]
where $U,V$ are constant unitary matrices, each $b_k(z)$ is a finite Blaschke product, and each $P_k$ is a projection matrix. Each factor $b_k(z)P_k+(I-P_k)$ belongs to $\rif(n,n)$. As noted above, since the unitary group is path connected we may assume $U=V=I_n$. Now let us write
\[
\begin{pmatrix} \Phi(z) \\ 0\end{pmatrix}  = \begin{pmatrix} b_1(z)P_1+(I-P_1) \\ 0 \end{pmatrix} \left(\prod_{k=2}^N (b_k(z)P_k +(I-P_k))\right).
\]
Let us work with
\begin{equation}\label{eqn:single-factor}
  \begin{pmatrix} b_1(z)P_1+(I-P_1) \\ 0 \end{pmatrix}.
\end{equation}
Conjugating by a unitary we may assume $b_1(z)P_1+(I-P_1)$ has the diagonal form
\[
\begin{pmatrix} b_1(z) & & & & & \\ & \ddots & & & & \\ & & b_1(z) & & & \\ & & & 1 & & \\ & & & & \ddots & \\ & & & & & 1\end{pmatrix}.
  \]
  Note that now, each column belongs to $\rif(n,1)$. Within $\rif(n+1,1)$ there is a path
  \[
  t\to \begin{pmatrix} (1-t)b_1(z) +t \\ 0 \\ \vdots \\ 0 \\ (\sqrt{t-t^2}) (1-b_1(z))\end{pmatrix}
  \]
  joining ${\begin{pmatrix} b_1(z) & 0 & \cdots & 0 & 0\end{pmatrix}}^T$ to ${\begin{pmatrix} 1 & 0 & \cdots & 0 & 0\end{pmatrix}}^T$. Doing this in the first column of the matrix (\ref{eqn:single-factor}) leaves the other columns unaffected and the whole path will lie in $\rif(m,n)$ (adding additional zeroes to the bottom of the column, if needed, to bring the size from $n+1$ up to $m$). We may thus successively move each diagonal entry $b_1(z)$ to $1$. Thus, our original $\begin{pmatrix} \Phi \\ 0\end{pmatrix}$ is now joined by a path in $\rif(m,n)$ to
    \[
    \begin{pmatrix} I_n \\ 0 \end{pmatrix} \left(\prod_{k=2}^N (b_k(z)P_k +(I-P_k))\right).
    \]
    We may then absorb the next Blaschke-Potapov factor into the column:
    \[
    \begin{pmatrix} b_2(z)P_2+(I-P_2) \\ 0 \end{pmatrix} \left(\prod_{k=3}^N (b_k(z)P_k +(I-P_k))\right)
    \]
    and repeat the process, so that in the end we see that $\begin{pmatrix} \Phi \\ 0\end{pmatrix}$ is joined to $\begin{pmatrix} I_n \\ 0\end{pmatrix}$ in $\rif(m,n)$ as desired. 
\end{proof}

\bibliographystyle{plain} 
\bibliography{homotopy-RIF}

\end{document}